\documentclass[A4paper,11pt]{article}
\usepackage{amsfonts}
\usepackage{amsmath}
\usepackage{amssymb}
\usepackage{graphicx}
\begin{document}
\title{Homological thickness and stability of torus knots  
\author{Marko Sto\v si\'c 
\thanks{The author was supported by {\it Funda\c c\~ao de Ci\^encia e
Tecnologia}/(FCT), grant no. SFRH/BD/6783/2001}\\
\\
Departamento de Matem\'atica \\
Instituto Superior T\'ecnico\\
Av. Rovisco Pais 1\\
1049-001 Lisbon\\ 
Portugal\\
e-mail: mstosic@math.ist.utl.pt
}
}

\date{}

\newtheorem{theorem}{Theorem}
\newtheorem{acknowledgment}[theorem]{Acknowledgment}
\newtheorem{algorithm}[theorem]{Algorithm}
\newtheorem{axiom}[theorem]{axiom}
\newtheorem{case}[theorem]{Case}
\newtheorem{claim}[theorem]{Claim}
\newtheorem{conclusion}[theorem]{Conclusion}
\newtheorem{condition}[theorem]{Condition}
\newtheorem{conjecture}[theorem]{Conjecture}
\newtheorem{corollary}[theorem]{Corollary}
\newtheorem{criterion}[theorem]{Criterion}
\newtheorem{definition}{Def\mbox{}inition}
\newtheorem{example}{Example}
\newtheorem{exercise}[theorem]{Exercise}
\newtheorem{lemma}{\indent Lemma}
\newtheorem{notation}[theorem]{Notation}
\newtheorem{problem}[theorem]{Problem}
\newtheorem{proposition}{Proposition}
\newtheorem{remark}[theorem]{Remark}
\newtheorem{solution}[theorem]{Solution}
\newtheorem{summary}[theorem]{Summary}
\newcommand{\ud}{\mathrm{d}}

\def\gcd{\mathop{\rm gcd}}
\def\Ker{\mathop{\rm Ker}}
\def\max{\mathop{\rm max}}
\def\map{\mathop{\rm map}}
\def\lcm{\mathop{\rm lcm}}
\def\kraj{\hfill\rule{6pt}{6pt}}
\def\diag{\mathop{\rm diag}}
\def\span{\mathop{\rm span}}
\def\deg{\mathop{\rm deg}}
\def\rank{\mathop{\rm rank}}
\def\sgn{\mathop{\rm sgn}}
\def\kvn{\{n\}_q}
\def\F{\mathbb{F}}
\def\R{\mathbb{R}}
\def\C{\mathcal{C}}
\def\P{\mathcal{P}}
\def\A{\mathcal{A}}
\def\D{\mathcal{D}}
\def\H{\mathcal{H}}
\def\N{\mathbb{N}}
\def\K{\mathbb{K}}
\def\Z{\mathbb{Z}}
\def\Q{\mathbb{Q}}
\def\X{\qbezier(0.00,0.00)(0.50,1.00)(1.00,2.00)
\qbezier(1.00,0.00)(0.80,0.40)(0.60,0.80)
\qbezier(0.00,2.00)(0.20,1.60)(0.40,1.20)
}

\def\Y{\qbezier(1.00,0.00)(0.50,1.00)(0.00,2.00)
\qbezier(0.00,0.00)(0.20,0.40)(0.40,0.80)
\qbezier(1.00,2.00)(0.80,1.60)(0.60,1.20)
}

\def\O{\qbezier(0.00,0.00)(0.20,1.00)(0.00,2.00)
\qbezier(1.00,0.00)(0.80,1.00)(1.00,2.00)
}
\def\OF{\qbezier(0.00,0.00)(0.00,1.00)(0.00,2.00)
\qbezier(1.00,0.00)(1.00,1.00)(1.00,2.00)
}

\def\l{
\qbezier(0.00,0.00)(0.50,1.20)(1.00,0.00)
\qbezier(0.00,2.00)(0.50,0.80)(1.00,2.00)
}

\def\LPP{
\qbezier(0.00,0.00)(0.35,0.70)(0.70,1.40)
\qbezier(0.70,0.60)(0.65,0.70)(0.60,0.80)
\qbezier(0.00,2.00)(0.20,1.60)(0.40,1.20)
\qbezier(0.7,1.4)(1.3,2.6)(1.3,1)
\qbezier(0.7,0.6)(1.3,-0.6)(1.3,1)}

\def\LPM{\qbezier(0.00,0.00)(0.2,0.40)(0.40,0.80)
\qbezier(0.60,1.20)(0.65,1.30)(0.70,1.40)
\qbezier(0.00,2.00)(0.35,1.30)(0.70,0.60)
\qbezier(0.7,1.4)(1.3,2.6)(1.3,1)
\qbezier(0.7,0.6)(1.3,-0.6)(1.3,1)}

\def\XS{\qbezier(0.00,0.00)(0.50,1.00)(1.00,2.00)
\qbezier(1.00,0.00)(0.80,0.40)(0.60,0.80)
\qbezier(0.00,2.00)(0.20,1.60)(0.40,1.20)

\qbezier(0.0,2.00)(0.0,1.85)(0,1.70)
\qbezier(1.00,2.00)(1,1.85)(1.0,1.70)
\qbezier(0.00,2.00)(0.10,1.90)(0.20,1.80)
\qbezier(1.00,2.00)(0.9,1.9)(0.8,1.80)}

\def\YS{\qbezier(1.00,0.00)(0.50,1.00)(0.00,2.00)
\qbezier(0.00,0.00)(0.20,0.40)(0.40,0.80)
\qbezier(1.00,2.00)(0.80,1.60)(0.60,1.2)
\qbezier(0.0,2.00)(0.0,1.85)(0,1.70)
\qbezier(1.00,2.00)(1,1.85)(1.0,1.70)
\qbezier(0.00,2.00)(0.10,1.90)(0.20,1.80)
\qbezier(1.00,2.00)(0.9,1.9)(0.8,1.80)}

\def\GR{\qbezier(0.00,0.00)(0.25,0.25)(0.50,0.50)
\qbezier(1.00,0.00)(0.75,0.25)(0.50,0.50)
\qbezier(0.00,2.00)(0.25,1.75)(0.50,1.5)
\qbezier(1.0,2.00)(0.75,1.75)(0.5,1.50)

\qbezier(0.00,2.00)(0,1.9)(0,1.80)
\qbezier(0.00,2.00)(0.1,2)(0.20,2)
\qbezier(1.00,2.00)(1,1.9)(1,1.80)
\qbezier(1.00,2.00)(0.9,2)(0.8,2)
\linethickness{2.5pt}
\qbezier(0.5,0.5)(0.5,1)(0.50,1.5)}

\def\GRF{\qbezier(0.00,0.00)(0.25,0.25)(0.50,0.50)
\qbezier(1.00,0.00)(0.75,0.25)(0.50,0.50)
\qbezier(0.00,2.00)(0.25,1.75)(0.50,1.5)
\qbezier(1.0,2.00)(0.75,1.75)(0.5,1.50)

\linethickness{2.5pt}
\qbezier(0.5,0.5)(0.5,1)(0.50,1.5)}

\def\OS{\qbezier(0.00,0.00)(0.60,1.00)(0.00,2.00)
\qbezier(1.00,0.00)(0.40,1.00)(1.00,2.00)
\qbezier(0.0,2.00)(0.0,1.85)(0,1.70)
\qbezier(1.00,2.00)(1,1.85)(1.0,1.70)
\qbezier(0.00,2.00)(0.10,1.90)(0.20,1.80)
\qbezier(1.00,2.00)(0.9,1.9)(0.8,1.80)}

\maketitle

\begin{abstract} 
In this paper we show that the non-alternating torus knots 
are homologically thick, i.e. that their Khovanov homology occupies at least three diagonals. Furthermore, we show that we can reduce the number of full twists of the torus knot without changing certain part of its homology, and consequently, we show that there exists stable homology of torus knots conjectured by Dunfield, Gukov and Rasmussen in \cite{dgr}. Since our main tool is the long exact sequence in homology, we have applied our approach in the case of the Khovanov-Rozansky ($sl(n)$) homology, and thus obtained analogous stability properties of $sl(n)$ homology of torus knots, also conjectured in \cite{dgr}.
\end{abstract}

\section{Introduction}

In recent years there has been a lot of interest in the 
 ``categorification" of link invariants, initiated by M. Khovanov in \cite{kov}. For each link $L$ in $S^3$ he def\mbox{}ined a graded chain complex, with grading preserving dif\mbox{}ferentials, whose graded Euler characteristic is equal to the Jones polynomial of the link $L$ (\cite{jones},\cite{kauf}), and whose homology groups (usually called $sl(2)$-homology groups) are link invariants. This is done by starting from the state-sum expression for the Jones polynomial (which is written as an alternating sum), then constructing for each term a graded module whose graded dimension is equal to the value of that term, and f\mbox{}inally, def\mbox{}ining the dif\mbox{}ferentials as appropriate grading preserving maps, so that the complex obtained is a link invariant (up to chain homotopy).

\indent Although the theory is rather new, it already has strong applications in low-dimensional topology. For instance, the short proof of the Milnor conjecture by Rasmussen in \cite{ras}, as well as the proof of the existence of exotic dif\mbox{}ferential structures on $\R^4$ (\cite{ras2}), which were previously accessible only by  
gauge theory.

\indent The advantage of Khovanov homology theory is that its def\mbox{}inition is combinatorial and since there is a straightforward algorithm for computing it, it is (theoretically) highly calculable. Nowadays there are several computer programs \cite{bnprog}, \cite{shum} that  can calculate ef\mbox{}fectively Khovanov homology of links with up to 50 crossings.

\indent Based on the calculations there are many conjectures about the properties of link homology, see e.g. \cite{bn}, \cite{pat}, \cite{dgr}. Some of the properties have been proved by now (see \cite{lee}, \cite{lee2}), but many of them are still open. \\

\indent In this paper we f\mbox{}irst show that the torus knots $T_{p,q}$ for $3\le p\le q$ (non-alternating torus knots) are homologically thick, i.e. that their Khovanov homology occupies at least three diagonals. Furthermore, in the course of the proof we obtain even stronger results that relate the homology of the torus knots $T_{p,q}$ and $T_{p,q+1}$. Namely, we prove that, up to a certain homological degree, their (unnormalized) homologies coincide. 

\indent As the f\mbox{}irst application of this result we calculate the homology of torus knots for low homological degrees. We also obtain the proof of the existence of stable Khovanov homology of torus knots, conjectured by Dunfield, Gukov and Rasmussen in \cite{dgr}.

\indent Furthermore, we conjecture that the homological width of the torus knot $T_{p,q}$ is at least $p$, and we reduce this problem to determining the nontriviality of certain homological groups.\\

\indent An analogous categorif\mbox{}ication of the $n$-specializations of the HOMFLYPT polynomial was carried out by M. Khovanov and L. Rozansky in 2004 (\cite{kovroz}). 
The construction uses the state-sum model for the HOMFLYPT polynomial (\cite{moy}) and is analogous to the categorif\mbox{}ication of the Jones polynomial: it uses the same cubic complex construction, and there exists a  similar long exact sequence in homology. However, since the state-sum model for the HOMFLYPT polynomial is much more complicated than  Kauf\mbox{}fman's state-sum model, the explicit calculation of the homology groups is practically impossible. Consequently, the values of the $sl(n)$-link homology are known only for a very small class of knots -- two-bridge knots (see \cite{ras2b}) and the closures of certain three-strand braids (\cite{web}).

\indent Since in the proofs of results for Khovanov homology of torus knots, we mainly use the long exact sequence of Khovanov homology (\ref{eq2}) -- we do not rely heavily on the explicit def\mbox{}inition of $sl(2)$ homology -- we also obtain most of the analogous results for the stability of $sl(n)$ homology  of torus knots.
Namely, we prove that up to a certain homological degree, the (unnormalized) $sl(n)$ homology groups of $T_{p,q}$ and $T_{p,q-1}$ torus knots coincide, and that there exists stable $sl(n)$-homology of torus knots, also conjectured in \cite{dgr}.\\

\indent The organization of the paper is as follows: in Section \ref{intro} we recall briefly the definition and the basic properties of Khovanov homology, we give short introduction to Khovanov-Rozansky homology and introduce notation for positive braid knots. In  Section \ref{glavni} we prove that the torus knots are homologically thick. In Section \ref{sledeci} we  relate the homologies of the torus knots $T_{p,q}$ and $T_{p,q+1}$, and we calculate the homology groups of torus knots in the homological degrees $0$, $1$, $2$, $3$ and $4$. Furthermore, in Section \ref{zadnji} we conjecture further results concerning the thickness of torus knots. Finally, in Section \ref{sln}, we prove that there exists stable $sl(n)$ homology of torus knots. \\

\noindent\textbf{Acknowledgements:} The author would like to thank M. Khovanov and J. Rasmussen for many helpful discussions and suggestions.

\section{Notation}\label{intro}

\subsection{Khovanov ($sl(2)$) homology}

We recall briefly the definition of Khovanov homology for links. For more details see \cite{bn}, \cite{kov}.  \\
\indent First of all, take a link $K$, its planar projection $D$, and take an ordering of the crossings of $D$. For each crossing $c$ of $D$, we define 0-resolution $D_0$ and 1-resolution $D_1$, as in the figure below. {{
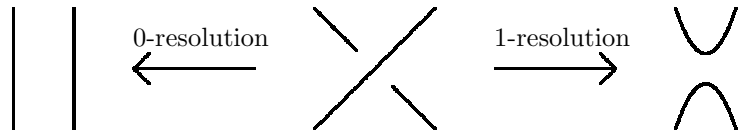
\begin{center}
\setlength{\unitlength}{4mm}
\begin{picture}(25,6)
\linethickness{0.8pt}
\qbezier(1.00,1.00)(1.00,3.00)(1.00,5.00)
\qbezier(3.00,1.00)(3.00,3.00)(3.00,5.00)

\qbezier(5.00,3.00)(7.00,3.00)(9.00,3.00)
\qbezier(5.00,3.00)(5.25,3.25)(5.50,3.50)
\qbezier(5.00,3.00)(5.25,2.75)(5.50,2.50)

\qbezier(15.00,5.00)(13.00,3.00)(11.00,1.00)
\qbezier(15.00,1.00)(14.30,1.70)(13.60,2.40)
\qbezier(11.00,5.00)(11.70,4.30)(12.40,3.60)

\qbezier(17.00,3.00)(19.00,3.00)(21.00,3.00)
\qbezier(21.00,3.00)(20.75,3.25)(20.50,3.50)
\qbezier(21.00,3.00)(20.75,2.75)(20.50,2.50)

\qbezier(23.00,5.00)(24.00,2.00)(25.00,5.00)
\qbezier(23.00,1.00)(24.00,4.00)(25.00,1.00)

\put(5.00,3.75){{\small\textrm{0-resolution}}}
\put(17.00,3.75){{\small\textrm{1-resolution}}}

\end{picture}
\end{center}
}}

Denote by $m$ the number of crossings of  $D$. Then there is bijective correspondence between the total resolutions of $D$ and the set $\{0,1\}^m$. Namely, to every $m$-tuple $\epsilon=(\epsilon_1,\ldots,\epsilon_m)\in \{0,1\}^m$ we associate the resolution $D_{\epsilon}$ where we resolved the $i$-th crossing in a $\epsilon_i$-resolution.\\
\indent Every resolution $D_{\epsilon}$ is a collection of disjoint circles. To each circle we associate graded $\Z$-module $V$, which is freely generated by two basis vectors $1$ and $X$, with $\deg 1 = 1$ and $\deg X=-1$. To $D_\epsilon$ we associate the module $M_{\epsilon}$, which is the tensor product of $V$'s over all circles in the resolution. Now, all the resolution $D_\epsilon$ with fixed $|\epsilon|$ (sum of elements of $\epsilon$) are grouped, and all resolutions are drawn as (skewed) $m$-dimensional cube such that in $i$-th column are the resolutions $D_{\epsilon}$ with $|\epsilon|=i$. The $i$-th chain group $C^i$ is given by:
$$C^i(D)=\oplus_{|\epsilon|=i}M_{\epsilon}\{i\}.$$
\indent Here, by $\{i\}$, we have denoted shift in grading of $M_{\epsilon}$ (for more details see e.g. \cite{bn}).\\
\indent The differential $d^i:C^i(G)\to C^{i+1}(G)$ is defined as (signed) sum of  ``per-edge" differentials. Namely, the only nonzero maps are from $D_{\epsilon}$ to $D_{\epsilon'}$, where $\epsilon=(\epsilon_1,\ldots,\epsilon_m)$, $\epsilon_i \in \{0,1\}$, if and only if $\epsilon'$ has all entries same as $\epsilon$ except one $\epsilon_j$, for some $j\in \{0,1\}$, which is changed from 0 to 1. We denote these differentials by $d_{\nu}$, where $\nu$ is $m$-tuple which consists of the label $\ast$ at the position $j$ and of $m-1$ 0's and 1's (the same as the remaining entries of $\epsilon$). 
Note that in these cases, either two circles of $D_{\epsilon}$ merge into one circle of $D_{\epsilon'}$ or one circle of $D_{\epsilon}$ splits into two circles of $D_{\epsilon'}$, and all other circles remain the same. \\
\indent In the first case, the map $d_{\nu}$ is defined as the identity on the tensor factors ($V$) that correspond to the unchanged circles, and on the remaining factors is given as the (graded preserving) multiplication map $m:V\otimes V \to V\{1\}$, which is given on basis vectors by:
$$m(1\otimes 1)=1,\quad m(1\otimes X)=m(X\otimes 1)=X,\quad m(X\otimes X)=0.$$
\indent In the second case, the map $d_{\nu}$ is defined as the identity on the tensor factors ($V$) that correspond to the unchanged circles,
and on the remaining factors is given as the (graded preserving) comultiplication map $\Delta:V\to V \otimes V\{1\}$, which is given on basis vectors by:
$$\Delta(1)=1\otimes X+X\otimes 1,\quad \Delta(X)=X\otimes X.$$
\indent Finally, to obtain the differential $d^i$ of the chain complex $C(D)$, we sum all contributions $d_{\nu}$ with $|\nu|=i$, multiplied by the sign $(-1)^{f(\nu)}$, where $f(\nu)$ is equal to the number of 1's ordered before $\ast$ in $\nu$. This makes every square of our cubic complex anticommutative, and hence we obtain the genuine differential (i.e. $(d^i)^2=0$).\\

\indent The homology groups of the obtained complex $(C(D),d)$ we denote by $H^i(D)$ and call \textit{unnormalized homology groups of $D$}. In order to obtain link invariants (i.e. independence of the chosen projection), we have to shift the chain complex (and hence the homology groups) by:
\begin{equation}
\C(D)=C(D)[-n_{-}]\{n_+-2n_-\},
\label{eq1}
\end{equation}
where $n_+$ and $n_-$ are the numbers of positive and negative crossings, respectively, of the diagram $D$ (see below for conventions). 
{{\begin{center}
\setlength{\unitlength}{4mm}
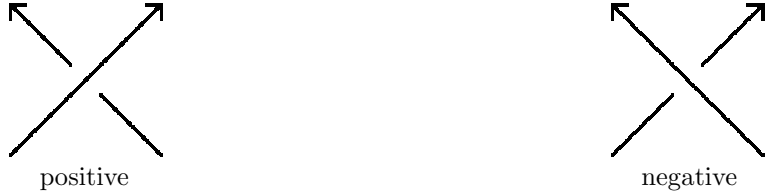
\begin{picture}(30,6) \label{sl1}
\linethickness{0.8pt}
\put(2.7,0){\qbezier(5.00,5.00)(2.50,2.50)(0.00,0.00)
\qbezier(5.00,0.00)(4.00,1.00)(3.00,2.00)
\qbezier(0.00,5.00)(1.00,4.00)(2.00,3.00)
\qbezier(0.00,5.00)(0.25, 5.00)(0.50,5.00)
\qbezier(0.00,5.00)(0.00,4.75)(0.00,4.50)
\qbezier(5.00,5.00)(4.75, 5.00)(4.50,5.00)
\qbezier(5.00,5.00)(5.00,4.75)(5.00,4.50)
\qbezier(22.00,2.00)(21.00,1.00)(20.00,0.00)
\qbezier(25.00,5.00)(24.00,4.00)(23.00,3.00)
\qbezier(20.00,5.00)(22.50,2.50)(25.00,0.00)
\qbezier(20.00,5.00)(20.25,5.00)(20.50,5.00)
\qbezier(20.00,5.00)(20.00,4.75)(20.00,4.50)
\qbezier(25.00,5.00)(24.75,5.00)(24.50,5.00)
\qbezier(25.00,5.00)(25.00,4.75)(25.00,4.50)
\put(1.00,-1.00){{\small\textrm{positive}}}
\put(21.00,-1.00){{\small\textrm{negative}}} }
\end{picture}
\end{center}
}}

\vskip 0.5cm
In the formula (\ref{eq1}), we have denoted by $[-n_-]$, the shift in homology degrees (again, for more details see \cite{bn}).\\
\indent The homology groups of the complex $\C(D)$ we denote by $\H^i(D)$. Hence, we have $\H^{i,j}(D)=H^{i+n_-,j-n_++2n_-}(D)$. 

\begin{theorem}(\cite{kov},\cite{bn})
The homology groups $\H(D)$ are independent of the choice of the planar projection $D$. Furthermore, the graded Euler characteristic of the complex $\C(D)$ is equal to Jones polynomial of the link $K$.
\end{theorem}
\indent Hence, we can write $\H(K)$, and we call $\H^i(K)$ \textit{the homology groups of the link $K$}.\\

\indent Let $D$ be a diagram of a link $L$ and let $c$ be one of its  crossings. Denote by $D_i$, $i=0,1$ the diagram that is obtained after performing an $i$-resolution of the crossing $c$. Then one can see that the complex $C(D)$ is in fact the mapping cone of a certain homomorphism $f:C(D_0)\to C(D_1)$ (which is basically given by maps $m$ and $\Delta$). Hence, there exists a  long exact sequence of (unnormalized) homology groups (see e.g. \cite{v}):

{\small{
\begin{equation}\cdots\rightarrow H^{i-1,j-1}(D_1)\to H^{i,j}(D)\to H^{i,j}(D_0)\to H^{i,j-1}(D_1)\to H^{i+1,j}(D)\to \cdots\label{eq2}\end{equation}
}}

This long exact sequence is the categorification of the defining recursive relation of the Kauffman bracket:
$$ \langle D \rangle = \langle D_0 \rangle - q\langle D_1 \rangle.$$
Indeed, this relation can be obtained by taking the graded Euler characteristic of (\ref{eq2}). \\

\indent If $K$ is a positive knot (the knot that has a planar  projection with only positive crossings) then $\H^i(K)$ is trivial for all $i<0$. Furthermore, if $D$ is planar projection of positive knot $K$, with $n_-$ negative crossings then $H^i(D)$ is trivial for $i<n_-$. \\

\indent Usually, the homology groups of the link $K$ are represented as a planar array in such a way that  $\rank \H^{i,j}(K)=\dim (\H^{i,j}(K)\otimes \Q)$ (or the whole group $\H^{i,j}(K)$, if we want to keep track of the torsions) is specif\mbox{}ied in the position $(i,j)$. As can be easily seen, the $q$-gradings ($j$) of the generators of nontrivial $\H^{i,j}(K)$ and the number of components of $K$ are of the same parity (either all are even or all are odd). Hence, by a \textit{diagonal} of the homology of the link $K$, we mean a line $j-2i=a=const$, when there exist integers $i$ and $j$ such that $j-2i=a$ and $\rank \H^{i,j}(K)>0$. If $a_{\max}$ and $a_{\min}$ are the maximal and minimal value of $a$ such that the line $j-2i=a$ is a diagonal of the homology of the link $K$, then we def\mbox{}ine the homological width of the link $K$ to be  $h(K)=(a_{\max}-a_{\min})/2+1$.\\
\indent Every knot (link) occupies at least two diagonals (i.e. $h(K) \ge 2$ for every link $K$), and the ones that occupy exactly two diagonals are called H-thin, or homologically thin. For example all alternating knots are H-thin (\cite{lee}), and the free part of the homology of any H-thin knot is determined by its Jones polynomial and the signature. A knot that is not H-thin is called H-thick or homologically thick.\\
\indent Furthermore, for an element $x\in \H^{i,j}(K)$, we denote its homological grading -- $i$ -- by $t(x)$, and its $q$-grading (also called the quantum grading) -- $j$ -- by $q(x)$. We also introduce a third grading $\delta(x)$ by $\delta(x)=q(x)-2t(x)$. Hence, we have that the knot $K$ is H-thick, if there exist three generators of $\H(K)$ with dif\mbox{}ferent values of the  $\delta$-grading.\\
\indent An alternative way of presenting the homology of the knot is by means of the two-variable Poincar\'e polynomial $P(K)(t,q)$ of the chain complex $\C(D)$, i.e.:
$$P(K)(t,q)=\sum_{i,j\in \Z}{t^iq^j \rank \H^{i,j}(K)}.$$

\subsection{Khovanov-Rozansky ($sl(n)$) homology}\label{kovroznot}

In \cite{kovroz} M. Khovanov and L. Rozansky generalized the construction from the previous subsection, to the case of $sl(n)$ specialization of HOMFLYPT polyomial, for every $n\in\N$. This is done by categorifying the Murakami-Ohtsuki-Yamada (MOY) calculus \cite{moy} -- the generalization of the Kauffman's state model for the Jones polynomial. The main form of the construction is the same: namely, they again assigned to every diagram $D$ the cubic complex of resolutions, to each vertex (total resolution) is assigned appropriate graded vector space and to every edge of the cube, certain graded preserving map. The first difference is that we start from the {\em oriented} diagram $D$, and 0- and 1-resolutions at a certain crossing $c$ are defined according to the sign of the crossing as on the following picture:
\vskip 0.4cm

{{
\setlength{\unitlength}{4mm}
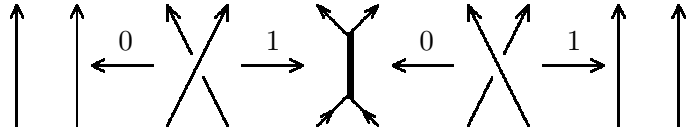
\begin{picture}(25,4.5) 
\linethickness{0.6pt}
\put(2,-1){
\qbezier(4.00,1.00)(4.00,3.00)(4.00,5.00)
\qbezier(3.80,4.60)(3.90,4.80)(4.00,5.00)
\qbezier(4.20,4.60)(4.10,4.80)(4.00,5.00) \qbezier(2.00,1.00)(2.00,3.00)(2.00,5.00)
\qbezier(1.80,4.60)(1.90,4.80)(2.00,5.00)
\qbezier(2.20,4.60)(2.10,4.80)(2.00,5.00) 
\qbezier(9.5,3)(10.5,3)(11.5,3)
\qbezier(11.5,3)(11.3,3.1)(11.1,3.2)
\qbezier(11.5,3)(11.3,2.9)(11.1,2.8)

\qbezier(4.5,3)(5.5,3)(6.5,3)
\qbezier(4.5,3)(4.7,3.1)(4.9,3.2)
\qbezier(4.5,3)(4.7,2.9)(4.9,2.8)

\put(10,0){
\qbezier(9,1)(8,3)(7,5)
\qbezier(7,1)(7.4,1.8)(7.8,2.6)
\qbezier(9,5)(8.6,4.2)(8.2,3.4)
}

\linethickness{2.5pt}
\qbezier(13.00,2.00)(13.00,3.00)(13.00,4.00)
\linethickness{0.6pt}
\qbezier(12.00,1.00)(12.50,1.50)(13.00,2.00)
\qbezier(14.00,1.00)(13.50,1.50)(13.00,2.00)
\qbezier(12.00,5.00)(12.50,4.50)(13.00,4.00)
\qbezier(14.00,5.00)(13.50,4.50)(13.00,4.00)

\qbezier(14.00,5.00)(13.80,4.90)(13.60,4.80)
\qbezier(14.00,5.00)(13.90,4.80)(13.80,4.60)

\qbezier(12.00,5.00)(12.20,4.90)(12.40,4.80)
\qbezier(12.00,5.00)(12.10,4.80)(12.20,4.60)

\qbezier(12.50,1.50)(12.30,1.40)(12.10,1.30)\qbezier(12.50,1.50)(12.40,1.30)(12.30,1.10)

\qbezier(13.50,1.50)(13.70,1.40)(13.90,1.30)\qbezier(13.50,1.50)(13.60,1.30)(13.70,1.10)
\qbezier(14.5,3)(15.5,3)(16.5,3)
\qbezier(14.5,3)(14.7,3.1)(14.9,3.2)
\qbezier(14.5,3)(14.7,2.9)(14.9,2.8)

\qbezier(19.5,3)(20.5,3)(21.5,3)
\qbezier(21.5,3)(21.3,3.1)(21.1,3.2)
\qbezier(21.5,3)(21.3,2.9)(21.1,2.8)
\put(-10,0){
\qbezier(19,5)(18,3)(17,1)
\qbezier(17,5)(17.4,4.2)(17.8,3.4)
\qbezier(19,1)(18.6,1.8)(18.2,2.6)

\qbezier(19,5)(19,4.75)(19,4.5)
\qbezier(19,5)(18.8,4.8)(18.6,4.6)
\qbezier(17,5)(17,4.75)(17,4.5)
\qbezier(17,5)(17.2,4.8)(17.4,4.6)
}
\put(10,0){
\qbezier(9,5)(9,4.75)(9,4.5)
\qbezier(9,5)(8.8,4.8)(8.6,4.6)
\qbezier(7,5)(7,4.75)(7,4.5)
\qbezier(7,5)(7.2,4.8)(7.4,4.6)
}

\qbezier(22.00,1.00)(22.00,3.00)(22.00,5.00)
\qbezier(21.80,4.60)(21.90,4.80)(22.00,5.00)
\qbezier(22.20,4.60)(22.10,4.80)(22.00,5.00)

\qbezier(24.00,1.00)(24.00,3.00)(24.00,5.00)
\qbezier(23.80,4.60)(23.90,4.80)(24.00,5.00)
\qbezier(24.20,4.60)(24.10,4.80)(24.00,5.00)

\put(5.40,3.50){0}
\put(10.3,3.50){1}
\put(15.40,3.50){0}
\put(20.3,3.50){1}
}
\end{picture}
\label{rezsln}
}}
\vskip 0.4cm

Hence, the total resolutions are in this case trivalent graph with thick edges (edges labelled 2 in \cite{moy}), and the values assigned to them in \cite{moy} satisfy certain set of (MOY) axioms. In \cite{kovroz}, the corresponding graded vector spaces are defined (in a rather complicated way) such  that they ``categorify" those axioms. Because of the complexity of this construction, the values of the $sl(n)$ homology are known only for very small set of knots.

\indent On the other hand, the main concepts and properties are the same as in the $sl(2)$ case - the cubic complex, mapping cone and consequently, the long exact sequence in $sl(n)$ homology: if $c$ is a positive crossing of an oriented diagram $D$, then there exists long exact sequence in (unnormalized) $sl(n)$ homology:

{\small{
\begin{equation}\cdots\rightarrow H_n^{i-1,j+1}(D_1)\to H^{i,j}_n(D)\to H^{i,j}_n(D_0)\to H^{i,j+1}_n(D_1)\to H^{i+1,j}_n(D)\to \cdots\end{equation}
}}
where $D_i$, $i=0,1$ is a diagram obtained from $D$ after resolving the crossing $c$ into an $i$-resolution. Obviously, in these long exact sequence we will always have diagrams of knots which also have trivalent vertices and thick edges, and the Khovanov-Rozansky prescription also assigns to them corresponding chain complexes and homology groups. We call such diagrams, the generalized regular diagrams. As in the $sl(2)$ case, the generalized regular diagrams with only positive crossings have trivial homology groups in negative homological degrees (since the chain groups in these degrees are trivial).

\indent For more details about the $sl(n)$ homology, we refer the reader to \cite{kovroz}.

\subsection{Positive braid knots}
  
The positive braid knots are the knots (or links) that are the closures of positive braids. Let $K$ be arbitrary positive braid knot and let $D$ be its planar projection which is the closure of a positive braid. Denote the number of strands of that braid by $p$. We say that  the crossing $c$ of $D$ is of the type $\sigma_i$, $i<p$, if it corresponds to the generator $\sigma_i$ in the braid word of which $D$ is the closure. Denote the number of crossings of the type $\sigma_i$ by $l_i$, $i=1,\ldots,p-1$ and order them from top to bottom. Then each crossing $c$ of $D$ we can write as the pair $(i,\alpha)$ (we will also write $(i\alpha)$ if there is no possibility of confusion), $i=1,\ldots,p-1$ and $\alpha=1,\ldots,l_i$, if $c$ is of the type $\sigma_i$ and it is ordered as $\alpha$-th among the crossings of the type $\sigma_i$. Finally, we order the crossings of $D$ by the following ordering: $c=(i\alpha)<d=(j\beta)$ if and only if $i<j$, or $i=j$ and $\alpha<\beta$.

\indent For some results on the homology of positive braid knots, see e.g. \cite{moj3}.

\section{Thickness of torus knots}\label{glavni}

\indent A knot or a link is a torus knot if it is isotopic to  a knot or a link that can be drawn without any points of intersection on the trivial torus. Every torus link is, up to a mirror image,  determined by two nonnegative integers $p$ and $q$, i.e. it is isotopic to a unique torus knot $T_{p,q}$ which has the diagram $D_{p,q}$ - the closure of the braid $(\sigma_1\sigma_2\ldots\sigma_{p-1})^q$ - as a planar projection. In other words, $D_{p,q}$ is the closure of the $p$-strand braid with $q$ full twists. Since $T_{p,q}$ is isotopic to $T_{q,p}$ we can assume that $p\le q$.

\indent  If $p=1$ then the torus knot $T_{p,q}$ is trivial and for $p=2$ the torus knot $T_{2,q}$ is alternating, hence its homology occupies exactly two diagonals. However, if $p\ge 3$, the torus knot $T_{p,q}$ is non-alternating and we will prove that its homology occupies at least three diagonals. Namely, we prove the following theorem:

\begin{theorem}\label{T1}
Let $K=T_{p,q}$, $3\le p\le q$ be a torus knot. Then 
$$\rank{\H^{4,(p-1)(q-1)+5}(K)}> 0.$$
\end{theorem}

From this theorem we obtain 
\begin{corollary}\label{C1}
Every torus knot $T_{p,q}$, $p,q\ge 3$ is H-thick, i.e. its Khovanov homology occupies at least three diagonals.
\end{corollary}
\textbf{Proof} (of Corollary \ref{C1}):\\

\indent Since $T_{p,q}$ is a positive knot, its zeroth homology group is two dimensional and the $q$-gradings (and consequently the $\delta$-gradings) of its two generators are $(p-1)(q-1)-1$ and $(p-1)(q-1)+1$, respectively (see e.g. \cite{ras}). However, from Theorem \ref{T1} we have that there exists a generator with $t$-grading equal to 4 and $q$-grading equal to $(p-1)(q-1)+5$, and so its $\delta$-grading is equal to $(p-1)(q-1)+5-2\cdot 4=(p-1)(q-1)-3$. Thus, we have obtained three generators    
of the homology of the torus knot $T_{p,q}$ which have three dif\mbox{}ferent values of the $\delta$-grading and hence its Khovanov homology occupies at least three diagonals. \kraj\\

\indent Now we give a proof of Theorem \ref{T1}.\\

\textbf{Proof}:\\

\indent 
First of all, since $K=T_{p,q}$ is a positive braid knot whose regular diagram $D_{p,q}$ is the closure of the braid $(\sigma_1\sigma_2\ldots\sigma_{p-1})^q$
with $(p-1)q$ crossings, we have that $\H^{4,(p-1)(q-1)+5}(K)=H^{4,6-p}(D_{p,q})$. So, we will ``concentrate" on calculating the latter homology group, i.e. showing that its rank is nonzero. In order to do this we will use the long exact sequence (\ref{eq2}) and we will relate the unnormalized fourth homology groups of the standard regular diagrams of the torus knots $D_{p,q}$ and $D_{p,q-1}$ for $p<q$.\\

\indent Let $3\le p <q$. Let $c_{p-1}$ be the crossing $(p-1,1)$ of the diagram $D_{p,q}$. Now denote by $E_{p,q}^1$ and $D_{p,q}^1$ the 1- and 0-resolutions, respectively, of the diagram $D_{p,q}$ at the crossing $c_{p-1}$. Then from (\ref{eq2}) we obtain the following long exact sequence
{\small{
\begin{equation*}\cdots\rightarrow H^{3,j-1}(E^1_{p,q})\to H^{4,j}(D_{p,q})\to H^{4,j}(D_{p,q}^1)\to H^{4,j-1}(E^1_{p,q})\to H^{5,j}(D_{p,q})\to \cdots\end{equation*}
}}
\indent Now, we can continue the process, and resolve the crossing $c_{p-2}=(p-2,1)$ of $D_{p,q}^1$ in two possible ways. Denote the diagram obtained by the 1-resolution by $E^2_{p,q}$, and the diagram obtained by the 0-resolution by $D^2_{p,q}$. Then from (\ref{eq2}) we have the long exact sequence
{\small{
\begin{equation*}\cdots\rightarrow H^{3,j-1}(E^2_{p,q})\to H^{4,j}(D^1_{p,q})\to H^{4,j}(D_{p,q}^2)\to H^{4,j-1}(E^2_{p,q})\to H^{5,j}(D^1_{p,q})\to \cdots\end{equation*}
}}
\indent After repeating this process $p-1$ times (resolving the crossing $c_{p-k}=(p-k,1),\,k=1,\ldots,p-1$, of $D^{k-1}_{p,q}$, obtaining the 1-resolution $E^k_{p,q}$ and 0-resolution $D^k_{p,q}$ and applying the same long exact sequence in homology), we obtain that for every $i=1,\ldots,p-1$, the following sequence is exact:
{\small{
\begin{equation}\cdots\rightarrow H^{3,j-1}(E^i_{p,q})\to H^{4,j}(D^{i-1}_{p,q})\to H^{4,j}(D_{p,q}^i)\to H^{4,j-1}(E^i_{p,q})\to H^{5,j}(D^{i-1}_{p,q})\to \cdots\label{eq3}\end{equation}
}}
Here  $D^0_{p,q}$ denotes $D_{p,q}$, and we obviously have that $D^{p-1}_{p,q}=D_{p,q-1}$.\\
\indent Our goal is to show $H^3(E^i_{p,q})$ and $H^4(E^i_{p,q})$ are trivial for every $0<i<p$. This is done in the following lemma.

\begin{lemma} \label{lem1}
For every three positive integers $p$, $q$ and $i$, such that $3 \le p < q$ and $i<p$, the knot with the diagram $E^i_{p,q}$ is positive, and the diagram $E^i_{p,q}$ has at least $p+q-3$ negative crossings.
\end{lemma}
\textbf{Proof:}\\

\indent Since for every $0<i<p$, $E^{i}_{p,q}$ is obtained by the 1-resolution of the crossing $c_{p-i}=(p-i,1)$ of the (positive braid knot) diagram $D^{i-1}_{p,q}$, it is the closure of the plat braid diagram with only one plat $E_{p-i}$:
\vskip 0.4cm
{\begin{center}\setlength{\unitlength}{4mm}
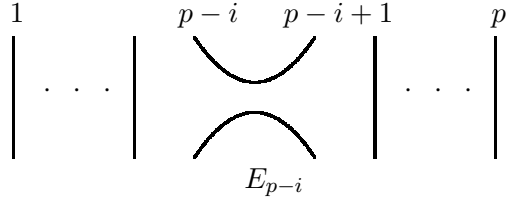
\begin{picture}(30,6)
\linethickness{0.8pt}
\put(1,0){
\qbezier(5.00,1.00)(5.00,3.00)(5.00,5.00)
\qbezier(9.00,1.00)(9.00,3.00)(9.00,5.00)
\qbezier(15.00,5.00)(13.00,2.00)(11.00,5.00)
\qbezier(11.00,1.00)(13.00,4.00)(15.00,1.00)
\qbezier(17.00,1.00)(17.00,3.00)(17.00,5.00)
\qbezier(21.00,1.00)(21.00,3.00)(21.00,5.00)
\put(6.00,3.00){$\cdot$}
\put(7.00,3.00){$\cdot$}
\put(8.00,3.00){$\cdot$}
\put(18.00,3.00){$\cdot$}
\put(19.00,3.00){$\cdot$}
\put(20.00,3.00){$\cdot$}
\put(4.90,5.50){1}
\put(10.50,5.50){$p-i$}
\put(14.00,5.50){$p-i+1$}
\put(20.90,5.50){$p$}
\put(12.70,0.00){$E_{p-i}$}
}
\end{picture}
\end{center}
}

 Now, note that two lower strands of $E_{p-i}$ are always ``neighbour" strands, i.e. they form a ribbon, through the diagram, until they reach the upper part of $E_{p-i}$. So, we can``slide" the lower part of the plat $E_{p-i}$ through the diagram (by using the second Reidemeister move -- R2 -- and also the f\mbox{}irst Reidemeister move -- R1 -- where the two strands intersect each other) until it reaches the left or the right hand side of the upper part of $E_{p-i}$. If it f\mbox{}irst reaches the right hand side, it automatically (or after a R1 move) becomes the closure of a positive braid diagram. If it f\mbox{}irst reaches the left hand side then after a R1 move and a ``slide" (sequence of R2 moves) we obviously obtain a positive braid diagram.\\
\indent Concerning the number of negative crossings of $E^i_{p,q}$, 
note that by performing the f\mbox{}irst sequence of R2 moves (sliding the lower part of $E_{p-i}$ through the diagram, from the top to the bottom) in each move we have ``canceled" one positive and one negative crossing. Furthermore, since $p<q$, the two strands of the lower part of the plat $E_{p-i}$ will make a full twist at least once, and so they will have two crossings with each other, which are both obviously negative crossings. So, we have that on each of the last (lower) $q-1$ blocks ($\sigma_1\ldots\sigma_{p-1}$) of $E^i_{p,q}$ we have at least one negative crossing. Furthermore, on the part where both lower strands of the plat $E_{p-i}$ make full twists, we have applied $p-2$ R2 moves, and hence we have in addition, at least, $p-2$ negative crossings.\\
\indent  Altogether, this gives at least $q-1+p-2=p+q-3$ negative crossings of $E^i_{p,q}$, as required. \kraj\\

\begin{remark}\label{rem1}
Obviously, the two ``neighbouring" strands from the Lemma above will make at least $[(q-1)/p]$ full twists, where by $[x]$ we have denoted the largest integer not greater than $x$. So, in fact we have proved that the diagram $E^i_{p,q}$ has at least $q-1+[(q-1)/p](p-2)$ negative crossings.
\end{remark}

\indent Now, we can go back to the proof. From the previous Lemma, we conclude that for every $0<k<p$, $\H^{i,j}(E^k_{p,q})$ is trivial for $i<0$, and that \begin{equation}
H^{i,j}(E^k_{p,q}) \textrm{ is trivial for } i<p+q-3,
\label{fspust}
\end{equation}
Hence, if $p+q-3>4$, from 
(\ref{eq3}) and (\ref{fspust}), we obtain that
$$H^{4,j}(D^{i-1}_{p,q})=H^{4,j}(D^i_{p,q}),\quad i=1,\ldots,p-1,$$
and thus we have that 
\begin{equation}
H^{4,j}(D_{p,q})=H^{4,j}(D_{p,q-1}),\quad\textrm {for } p+q>7, \,\, p<q.
\label{spust}
\end{equation}
So, we can decrease the number of full twists, $q$, without changing the fourth homology group.\\
\indent If $p=q=3$, then e.g. by using programs for computing Khovanov homology (\cite{bnprog}, \cite{shum}), we obtain that $\rank H^{4,3}(D_{3,3})=\rank \H^{4,9}(T_{3,3})=1$, as wanted. \\
\indent If $3=p<q$, then from (\ref{spust}), we have that $H^{4,3}(D_{3,q})=H^{4,3}(D_{3,4})$. However, by using programs for computing Khovanov homology, we obtain that the rank of the latter group (which is equal to $\H^{4,11}(T_{3,4})$) is equal to 1. 
\begin{remark}
In Bar-Natan's tables of knots in \cite{bn}, the link $T_{3,3}$ is denoted by $6^3_3$, and the torus knot $T_{3,4}$ is isotopic to the knot $8_{19}$. For the general notation of knots and links see \cite{rol} and \cite{bnprog}.
\end{remark}
\indent Now, let us move to the general case $4\le p\le q$. Then if $p<q$ we can apply (\ref{spust}) and obtain $H^{4,6-p}(D_{p,q})=H^{4,6-p}(D_{p,p})$. Thus, we are left with proving that the latter homology group is of nonzero rank.\\
\indent Now, apply the set of long exact sequences ($\ref{eq3}$) for the case $p=q$. In this case, like in Lemma \ref{lem1}, we obtain that $E^k_{p,p}$ is the diagram of a positive knot, for every $k=1,\ldots,p-1$. Furthermore, every diagram $E^k_{p,p}$ has exactly $2p-3$ negative crossings, and so we have:
\begin{equation}
H^i(E^k_{p,p}) \textrm{ is trivial for every }i<2p-3.
\label{novspust}
\end{equation} Since $p>3$, we have that $H^3(E^k_{p,p})$ and 
$H^4(E^k_{p,p})$ are trivial for all $k$, and so we have that
\begin{equation}
H^4(D_{p,p})=H^4(D_{p,p-1}).
\label{sp1}
\end{equation}
\indent On the other hand, we have that 
$$H^{4,6-p}(D_{p,p-1})=\H^{4,p^2-3p+7}(T_{p,p-1})=\H^{4,p^2-3p+7}(T_{p-1,p})=H^{4,7-p}(D_{p-1,p}).$$
\indent If $p>4$, then we have from (\ref{spust}) that $H^{4,7-p}(D_{p-1,p})=H^{4,7-p}(D_{p-1,p-1})$. 
By repeating this process, we can decrease the number of strands $p$, and obtain that:
$$H^{4,6-p}(D_{p,p})=H^{4,2}(D_{4,4}).$$
Finally, from (\ref{sp1}) we have $$H^{4,2}(D_{4,4})=H^{4,2}(D_{4,3})=\H^{4,11}(T_{3,4}),$$
and the last homology group, as we saw previously, is of rank 1. This concludes our proof.   \kraj 

\section{Stability of Khovanov homology for torus knots}\label{sledeci}

In the course of proving Theorem \ref{T1}, apart from showing that the torus knots are H-thick, we have obtained some other properties of the homology of the torus knots. Namely, we proved that we can reduce the number of full twists, $q$, of the standard diagram $D_{p,q}$ of the torus knot $T_{p,q}$ without changing the f\mbox{}irst $p+q-3$ homology groups. In other words, we have obtained the existence of stable homology of torus knots (see Section \ref{sln} below and \cite{dgr}).

\indent In the following theorem we summarize the stability properties obtained in the previous section.
\begin{theorem} \label{T3}
Let $p$, $q$ and $i$ be integers such that $2\le p<q$ and $i<p+q-3$. Then for every $j\in\Z$
\begin{equation}
H^{i,j}(D_{p,q})=H^{i,j}(D_{p,q-1}).
\label{f1}
\end{equation}
Furthermore, for every $2\le p<q$ and $i<2p-1$ and $j\in\Z$ we have 
\begin{equation}
H^{i,j}(D_{p,p+1})=H^{i,j}(D_{p,p+2})=\cdots=H^{i,j}(D_{p,q}).
\label{f2}
\end{equation}
Also, for every $p\ge 2$, $i<2p-3$ and $j\in\Z$, we have
\begin{equation}
H^{i,j}(D_{p,p})=H^{i,j+1}(D_{p-1,p}).
\label{f3}
\end{equation}
\end{theorem} 
\textbf{Proof:}\\
\indent The equations (\ref{f1}) and (\ref{f3}) we have already obtained in the course of proving Theorem \ref{T1} (long exact sequences (\ref{eq3}) with the homological degrees $3$, $4$ and $5$ replaced by $i-1$, $i$ and $i+1$ for every $i\in\Z$, respectively, and formulas (\ref{fspust}) and (\ref{novspust})). Formula (\ref{f2}) obviously follows from (\ref{f1}) since $i<2p-1=p+(p+2)-3$. \kraj\\

\begin{remark}
Bearing in mind Remark \ref{rem1}, in fact we have obtained that if $p$ and$q$ are integers such that $2\le p<q$ then for every $j\in\Z$ we have
\begin{equation}
H^{i,j}(D_{p,q})=H^{i,j}(D_{p,q-1}),\quad \mathrm{ for }\quad i<q-1+[(q-1)/p](p-2),
\end{equation}
\label{rem2}
\end{remark}

The torus knots $T_{2,q}$ are alternating and their homology is well-known (see e.g. \cite{kov}). However, as the f\mbox{}irst corollary of the previous theorem we obtain the homology groups of $T_{p,q}$, for $3\le p \le q$, with low homological degree. 

\begin{theorem}
Let $3\le p\le q$ with $p$ and $q$ not both equal to 3. Then we have 
\begin{eqnarray*}
\H^{0,(p-1)(q-1)\pm 1}(T_{p,q})&=&\Z\\
\H^{2,(p-1)(q-1)+3}(T_{p,q})&=&\Z\\
\H^{3,(p-1)(q-1)+7}(T_{p,q})&=&\Z\\
\H^{3,(p-1)(q-1)+5}(T_{p,q})&=&\Z_2\\
\H^{4,(p-1)(q-1)+6 \pm 1}(T_{p,q})&=&\Z.
\end{eqnarray*}
All other $\H^{i,j}(T_{p,q})$ for $i=0,\ldots,4$, are trivial. 
\end{theorem}
\textbf{Proof:}\\
\indent  Suppose that $p=3$. Then by applying  (\ref{f2}), we obtain that $H^{i,j}(D_{3,q})=H^{i,j}(D_{3,4})$ for $i=0,\ldots,4$.
If $p>3$, then by applying  (\ref{f1}) repeatedly, we obtain
$H^{i,j}(D_{p,q})=H^{i,j}(D_{p,p})$.
Furthermore, by applying (\ref{f3}) (and then (\ref{f1})) repeatedly we obtain $H^{i,j}(D_{p,p})=H^{i,j+p-3}(D_{3,4})$ for $i=0,1,2,3,4$.
Finally, the homology of the last torus knot is well-known, see e.g. \cite{shtor}-knot $8_{19}$, and thus we obtain the required result. \kraj\\

\section{Further thickness results}\label{zadnji}

Even though we have shown that the torus knots $T_{p,q}$, $p\ge 3$ are H-thick, from the existing experimental results one can see that the homology of torus $T_{p,q}$ knots occupies at least $p$ diagonals (i.e. that its homological width is at least $p$). In fact, one can see that in all  examples we have that $H^{2p-2,p}(D_{p,q})$ is of nonzero rank. 
\begin{proposition}
If $\rank{H^{2p-2,p}(D_{p,q})}>0$ then the homological width of the torus knot $T_{p,q}$ is at least $p$.
\end{proposition}
\textbf{Proof:}\\

\indent As we know, see e.g. the proof of Corollary \ref{C1}, there exists a generator of the homology group $\H^{0,(p-1)(q-1)+1}$ and its $\delta$-grading is equal to $(p-1)(q-1)+1$. Since we have assumed that   
$$\rank{H^{2p-2,p}(D_{p,q})}=\rank{\H^{2p-2,p+(p-1)q}(T_{p,q})}>0,$$
we have that there exists a generator of this homology group whose $\delta$-grading is equal to $p+(p-1)q-2(2p-2)=(p-1)(q-1)+3-2p$. So, we have two generators whose $\delta$-gradings dif\mbox{}fer by $2p-2$, and hence they lie on two dif\mbox{}ferent diagonals between which there are $p-2$ diagonals. Hence the homological width of the torus knot $T_{p,q}$ is at least $p$. \kraj\\

Thus we are left with proving that $\rank{H^{2p-2,p}(D_{p,q})}>0$ . From (\ref{f2}) we have that ${H^{2p-2,p}(D_{p,q})}={H^{2p-2,p}(D_{p,p+1})}$. Furthermore we have
\begin{lemma} 
${H^{2p-2,p}(D_{p,p})}={H^{2p-2,p}(D_{p,p+1})}.$
\end{lemma} 
\textbf{Proof:}\\

\indent In order to prove this, we will start from the diagram $D_{p,p+1}$ and we will use the same process as in the proof of Theorem \ref{T1}. Namely, we obtain the long exact sequences, see (\ref{eq3}), for every
$i=1,\ldots,p-1$:
\arraycolsep 0pt
{{
\begin{eqnarray}\cdots\rightarrow H^{2p-3,p-1}(E^i_{p,p+1})\to  H^{2p-2,p}(D^{i-1}_{p,p+1})\to \nonumber \\ \to  H^{2p-2,p}(D_{p,p+1}^i) \to H^{2p-2,p-1}(E^i_{p,p+1})\to \cdots\label{eq4}\end{eqnarray}
}}
\arraycolsep 6pt
For every $i=1,\ldots,p-1$ we can calculate explicitly the number of positive and negative crossings of $E^i_{p,p+1}$, and we can f\mbox{}ind explicitly the positive diagram to which $E^i_{p,p+1}$
is isotopic.\\

\indent One can easily see that the number of negative crossings of each $E^i_{p,p+1}$ is equal to $2p-2$ and hence the number of positive crossings is equal to $(p-1)(p+1)-i-(2p-2)=p^2-2p+1-i$. On the other hand, every $E^i_{p,p+1}$ for $i=1,\ldots,p-2$ is isotopic (by a sequence of R2 and R1 moves as explained in the proof of Lemma \ref{lem1}) to  the diagram  $D_{p-2,p-1}^{i-1}$, while $E^{p-1}_{p,p+1}$ is isotopic to the diagram $D_{p-2,p-1}^{p-3}\coprod U=D_{p-2,p-2}\coprod U$, where by $U$ we denote the unknot.  Hence we have that $\H^{l}(E_{p,p+1}^i)=0$ for $l<0$. Also, since $D_{p,q}^i$ is a positive braid knot with $p$ strands and $(p-1)q-i$ crossings, we have that  $\H^{0,(p-3)(p-2)-(i-1)\pm 1}(D_{p-2,p-1}^{i-1})=\Z$ and all other $\H^{0,j}(D_{p-2,p-1}^{i-1})$ are trivial. Hence, we have that the only nontrivial part of the zeroth homology group  of $E_{p,p+1}^i$ is given by 
 $\H^{0,(p-3)(p-2)-(i-1)\pm 1}(E_{p,p+1}^{i})=\Z$ for $i=1,\ldots,p-2$, and $\H^{0,(p-3)(p-2)-(p-3)\pm 1\pm 1}(E_{p,p+1}^{p-1})=\Z$. Thus,  we have that for every $i=1,\ldots,p-1$,  
$\H^{0,p^2-5p+4-i}(E_{p,p+1}^i)$ is trivial.\\

\indent Finally, since 
the number of negative crossings of $E_{p,p+1}^i$ is equal to $2p-2$, we have that $H^{2p-3}(E_{p,p+1}^i)$ is trivial. Furthermore, since the number of positive crossings  of $E_{p,p+1}^i$ is equal to $p^2-2p+1-i$ we have that 
$$H^{2p-2,p-1}(E_{p,p+1}^i)=\H^{0,p^2-5p+4-i}(E_{p,p+1}^i)$$
which is trivial for every $i=1,\ldots,p-1$. \\
\indent 
Hence from the long exact sequences (\ref{eq4}) we obtain 
$${H^{2p-2,p}(D_{p,p})}={H^{2p-2,p}(D_{p,p+1})},$$ as required. \kraj\\

\begin{conjecture} \label{con1}
The rank of the homology group $H^{2p-2,p}(D_{p,p})$ (and equivalently of $H^{2p-2,p}(D_{p,p+1})$) is nonzero.
\end{conjecture}

\indent As we saw, the validity of  Conjecture \ref{con1} implies that the homological width of the torus knot $T_{p,q}$ is at least $p$.\\
\indent Even though we don't (yet) have the proof of  Conjecture \ref{con1}, there is evidence that it is true. First of all, the computer program calculations show that the conjecture is true at least for $p \le 7$ (the calculations are mainly for knots, i.e. for $D_{p,p+1}$). Furthermore, Lee's variant $H_L^{i,j}$ of Khovanov homology (\cite{lee2}) for the $p$-component link $D_{p,p}$ has $2p$ generators in the homological degree $2p-2$. Also, as it is well-known, there exist spectral sequences whose $E_{\infty}$-page is Lee's homology  and whose $E_2$-page is Khovanov homology (see \cite{ras}, \cite{turner}). So $H_L^{i,j}\subset H^{i,j}$ and hence $H^{2p-2}(D_{p,p})$ has at least $2p$ generators. So, we are left with proving that at least one of them has the $q$-grading equal to $p$.\\

\section{Stable $sl(n)$ homology of torus knots} \label{sln}

\indent Def\mbox{}ine the following normalization of the Poincar\'e polynomial of the homology of the torus knot:
$$P_{m,n}(t,q)=q^{-(m-1)n}P(T_{m,n})(t,q).$$ 
Then from the ``descending" properties of Theorem \ref{T3} we have the following:
\begin{theorem} For every $m \in \N$ there exists a stable homology polynomial $P^S_m$ given by:
$$P^S_m(t,q)=lim_{n\to \infty}P_{m,n}(t,q).$$
\label{stablet}
\end{theorem}

\indent Furthermore, as we have shown, the (normalized) Poincar\'e polynomial $P_{m,n}(t,q)$ of the torus knot $T_{m,n}$ coincides with the stable polynomial $P^S_m$, for all powers of $t$ up to $m+n-3$.\\
\indent Similar results are obtained at the conjectural level in \cite{dgr} (with a conjectural bound on the powers of $q$ for  agreement between the stable homology and the ef\mbox{}fective homology of any  particular torus knot). In \cite{dgr},  reduced homology (see e.g. \cite{pat}) is used, but the whole method and all proofs from the previous sections work in the same way for reduced homology. \\

\indent Also, in \cite{dgr} the existence of stable $sl(n)$ homology for torus knots is conjectured. However, in the course of proving the stability property in the $sl(2)$ case (Theorem \ref{T3}, formula (\ref{f1})) the basic ingredient is the long exact sequence in homology together with the form of cube of resolutions. Since the analogous long exact sequence exists for $sl(n)$ homology (it is again the mapping cone), we can repeat the major part of the process. The long exact sequence in the case of $sl(n)$ homology is:

{\small{
\begin{equation}\cdots\rightarrow H^{i-1,j+1}_n(D_1)\to H^{i,j}_n(D)\to H^{i,j}_n(D_0)\to H^{i,j+1}_n(D_1)\to H^{i+1,j}_n(D)\to \cdots\label{lessln}\end{equation}
}}
where $D_i$, $i=0,1$ is obtained from $D$ after resolving the positive crossing $c$ into an $i$-resolution. Note that one of the diagrams $D_i$ is not a planar projection of a knot since it contains one thick edge. In the case that we are interested in (torus knots - positive knots), the diagram $D_1$ is the one which has one thick edge (for the details and notation see Section \ref{kovroznot} and \cite{kovroz}).

\indent Again, we start from the diagram $D_{p,q}$ of the torus knot $T_{p,q}$, and we resolve the crossing $c_{p-1}=(p-1,1)$. We denote the diagram obtained by the 0-resolution by $D^1_{p,q}$, and the diagram obtained by the $1$-resolution by $\bar{E}^1_{p,q}$. Then we have the following long exact sequence:

{\small{
\begin{equation*}\cdots\rightarrow H_n^{i-1,j+1}(\bar{E}^1_{p,q})\to H^{i,j}_n(D_{p,q})\to H_n^{i,j}(D_{p,q}^1)\to H_n^{i,j+1}(\bar{E}^1_{p,q})\to \cdots\end{equation*}
}}
We continue the process, by resolving the crossings $c_l=(l,1)$, $l=p-2,\ldots,1$ of the diagram $D^{p-1-l}_{p,q}$ and we denote the 0- and 1-resolution obtained, by $D^{p-l}_{p,q}$ and $\bar{E}_{p,q}^{p-l}$, respectively. Then we have the following long exact sequence

\begin{eqnarray*}
\cdots\rightarrow H_n^{i-1,j+1}(\bar{E}^l_{p,q})\to H^{i,j}_n(D^{l-1}_{p,q})\to H_n^{i,j}(D^l_{p,q})\to H_n^{i,j+1}(\bar{E}^l_{p,q})\to \cdots,\\
l=2,\ldots,p-1.
\end{eqnarray*}


Like in the $sl(2)$ case, we shall prove the following
\begin{lemma}\label{lem2} The homology group $H_n^i(\bar{E}^l_{p,q})$ is trivial for every $l<p$ and $i<p+q-3$.
\end{lemma}
This lemma, together with the above long exact sequences and the fact that $D_{p,q}^{p-1}=D_{p,q-1}$ gives
\begin{equation}
\label{spustsln}
H_n^{i,j}(D_{p,q-1})=H_n^{i,j}(D_{p,q}),\quad i<p+q-3.
\end{equation}
From this formula, we conclude 
the existence of the limit:
\begin{eqnarray*}
P^n_k(t,q)&=&\lim_{l\to\infty} \sum_{i,j\in \Z}t^iq^j \dim H_n^{i,j}(D_{k,l})=\\
&=&\lim_{l\to\infty} \sum_{i,j\in \Z}t^iq^j q^{(n-1)(k-1)l}\dim \H_n^{i,j}(T_{k,l})=\\
&=&\lim_{l\to\infty}q^{(n-1)(k-1)l}P^n(T_{k,l})(t,q),
\end{eqnarray*}
for every $k$, where $P^n(T_{k,l})(t,q)$ is the Poincar\'e polynomial of the chain complex assigned to $T_{k,l}$ by $sl(n)$-homology. In other words, we obtain
\begin{theorem}
There exists stable $sl(n)$ homology for torus knots.
\label{torussln}
\end{theorem}

\indent Thus, we are left with proving Lemma \ref{lem2}. 
We will use more or less the same approach as in Lemma \ref{lem1}. Let $C_n(\bar{E}^l_{p,q})$ be the chain complex assigned by $sl(n)$-link homology (\cite{kovroz}) to  $\bar{E}^l_{p,q}$. Then $H_n^i(\bar{E}^l_{p,q})=H^i(C_n(\bar{E}^l_{p,q}))$. Note that in the $sl(n)$ case, a complex $C_n$ is assigned to  generalized regular diagrams, i.e. to  regular diagrams where we also allow trivalent vertices and thick edges.

\indent Since $\bar{E}^l_{p,q}$ has only positive crossings, its homology groups are trivial for negative homological degrees. However, we will show that 
\begin{equation}
C_n(\bar{E}^l_{p,q})\sim D_n[p+q-3],
\label{slnkompl}
\end{equation}
where $D_n$ is a complex such that all its chain groups $D_n^i$ are trivial for $i<0$ (in fact, we will def\mbox{}ine $D_n$ as the direct sum of the complexes of the form $C_n(D_{\Gamma}^i)$, where $D_{\Gamma}^i$'s are the generalized regular diagrams whose all crossings are positive). Here, by $\sim$ we denote a quasi-isomorphism, which implies that the two complexes have isomorphic homology groups. Then (\ref{slnkompl}) implies that $H_n^i(\bar{E}^l_{p,q})=H^i(C_n(\bar{E}^l_{p,q}))$ is trivial for $i<p+q-3$.\\
\indent Like in Lemma \ref{lem1}, we have that the two lower strands (thin edges) of the $\bar{E}_i$ part, 
will form at least two crossings with each other
(corresponding to an R1 move in the proof of Lemma \ref{lem1}) and both of them will have at least $p+q-5$ over- or undercrossings with the same strand (corresponding to an R2 move in the proof of Lemma \ref{lem1}). We will show that in each of these cases we can ``shift" up our complex by one homological degree, and thus obtain (\ref{slnkompl}). \\
\indent In order to prove this we will use the following fact (``cancellation principle" for chain complexes): if we quotient the chain complex $\C$ by an (arbitrary) acyclic subcomplex $\C'$ (i.e. a subcomplex with trivial homology), then the  quotient complex $\C/\C'$ is quasi-isomorphic to the complex $\C$, and so they have isomorphic homology groups (see e.g. Lemma 3.7 of \cite{bn}).\\

\paragraph{Untwisting an R1 move}

\indent First, let us work with the analog of the R1 move. Let $\bar{D}$ be the diagram that contains the following diagram as a subdiagram:

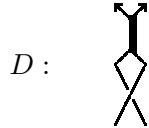
\begin{center}
\setlength{\unitlength}{4mm}
\begin{picture}(20,3.4) 
\linethickness{0.6pt}
\put(10,-0.2){
\put(0,0){\X}
\put(0,2){\GR}
}
\put(6.5,1.5){$D:$}
\end{picture}
\end{center}
Then (see \cite{kovroz}) the chain complex $C_n(\bar{D})$ associated to $\bar{D}$ is the mapping cone of a certain homomorphism $f:C_n(\bar{D}_0)\to C_n(\bar{D}_1)\{-1\}$, where $\bar{D}_0$ and $\bar{D}_1$ are the 0- and 1-resolutions, respectively, of the crossing of $D$, i.e. they look the same as the diagram $\bar{D}$ except that its subdiagram $D$ is replaced by $D_0$ and $D_1$, respectively:

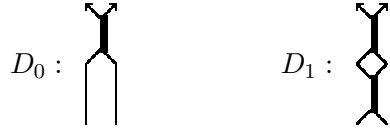
\begin{center}
\setlength{\unitlength}{4mm}
\begin{picture}(20,3.8) 
\linethickness{0.6pt}
\put(5,-0.2){
\put(0,0){\OF}
\put(0,2){\GR}
}
\put(14,-0.2){
\put(0,0){\GRF}
\put(0,2){\GR}
}
\put(2.5,1.5){$D_0:$}
\put(11.5,1.5){$D_1:$}
\end{picture}
\end{center}
In other words, the complex associated to $\bar{D}$ is the total complex of the complex given by $C_n(\bar{D}_0)\to C_n(\bar{D}_1)[1]\{-1\}$. Furthermore, we have that (Proposition 30 in \cite{kovroz} -- ``categorification" of one of the MOY axioms):
$$C_n(\bar{D}_1)\cong C_n(\bar{D}_0)\{1\} \oplus C_n(\bar{D}_0)\{-1\}.$$
We also have that the projection of $f$ to the f\mbox{}irst summand is an isomorphism (see \cite{kovroz} - invariance under the Reidemeister IIa move), and hence the complex $C_n(\bar{D})$ is quasi-isomorphic to $C_n(\bar{D}_0)[1]\{-2\}$ (by the cancellation principle). So, the last two complexes have isomorphic homology groups. Thus, we can ``untwist" the crossing involving two strands that are connected to the same thick edge (the analog of an R1 move in the $sl(2)$ case) by shifting the complex of the diagram obtained up by one in homological degree, as required. \\

\paragraph{Untwisting an R2 move}

\indent Hence, after untwisting  the two crossings of $\bar{E}^l_{p,q}$ that were resolved in the $sl(2)$ case by the R1 move, we are left with a diagram of the form:

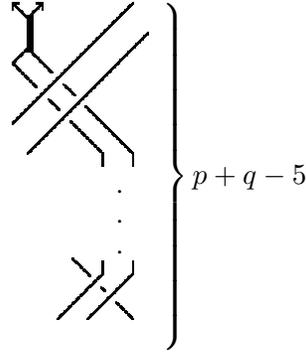
\begin{equation*}
\setlength{\unitlength}{4mm}
\left.
{
\begin{picture}(5,6) 
\linethickness{0.6pt}
\put(0,0){
\put(0,4){\GR}
\qbezier(0,2)(2,4)(4,6)
\qbezier(0.5,1)(2.5,3)(4.5,5)
\qbezier(0,4)(0.4,3.6)(0.8,3.2)
\qbezier(1,4)(1.15,3.85)(1.3,3.7)
\qbezier(1.2,2.8)(1.38,2.62)(1.56,2.44)
\qbezier(1.7,3.3)(1.88,3.12)(2.06,2.94)
\qbezier(1.94,2.06)(2.47,1.53)(3,1)
\qbezier(2.44,2.56)(3.22,1.78)(4,1)
\qbezier(3.44,-3.94)(3.72,-4.22)(4,-4.5)
\qbezier(3.06,-3.56)(3,-3.5)(2.94,-3.44)
\qbezier(2,-2.5)(2.28,-2.78)(2.56,-3.06)
\qbezier(3,-3)(2.25,-3.75)(1.5,-4.5)
\qbezier(4,-3)(3.25,-3.75)(2.5,-4.5)
\qbezier(3,-3)(3,-2.8)(3,-2.6)
\qbezier(4,-3)(4,-2.8)(4,-2.6)
\put(0,1){
\qbezier(3,-0.4)(3,-0.2)(3,0)
\qbezier(4,-0.4)(4,-0.2)(4,0)
}
}
\put(3.45,-0.5){$\cdot$}
\put(3.45,-1.5){$\cdot$}
\put(3.45,-2.5){$\cdot$}
\end{picture}}
\right\}
p+q-5
\end{equation*}


In other words, we have two neighbouring strands that are both connected to the same thick edge, and both go over or both go under $p+q-5$ strands. We will show that the complex corresponding to this diagram is quasi-isomorphic to the complex of diagrams whose crossings are all positive, shifted up in homological degree by $p+q-5$.\\
\indent Let $\bar{D}$ be a positive diagram that contains the following diagram as a subdiagram:

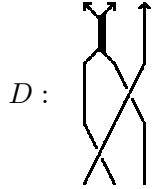
\begin{center}
\setlength{\unitlength}{4mm}
\begin{picture}(20,5.4) 
\linethickness{0.6pt}
\put(9,-0.2){
\put(0,0){\X}
\put(0,4){\GR}
\put(1,2){\X}
\qbezier(0,2)(0,3)(0,4)
\qbezier(2,6)(2,5)(2,4)
\qbezier(2,2)(2,1)(2,0)
\put(2,6){
\qbezier(0,0)(0.08,-0.08)(0.16,-0.16)
\qbezier(0,0)(-0.08,-0.08)(-0.16,-0.16)
}
}
\put(6.5,2.5){$D:$}
\end{picture}
\end{center}

\indent Denote by $D_0$, $D_1^1$, $D_1^2$ and $D_2$ the resolutions obtained from $D$ by resolving its two crossings, according to the following pictures:

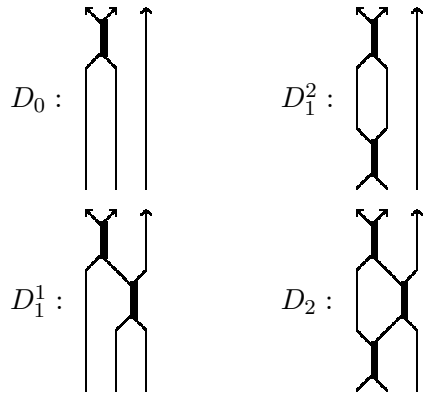
\begin{center}
\setlength{\unitlength}{4mm}
\begin{picture}(20,5.8) 
\linethickness{0.6pt}
\put(5,-0.2){
\put(0,0){\OF}
\put(0,2){\OF}
\put(0,4){\GR}
\qbezier(2,6)(2,3)(2,0)
\put(2,6){
\qbezier(0,0)(0.08,-0.08)(0.16,-0.16)
\qbezier(0,0)(-0.08,-0.08)(-0.16,-0.16)
}
}
\put(14,-0.2){
\put(0,0){\GRF}
\put(0,4){\GR}
\put(0,2){\OF}
\qbezier(2,6)(2,3)(2,0)
\put(2,6){
\qbezier(0,0)(0.08,-0.08)(0.16,-0.16)
\qbezier(0,0)(-0.08,-0.08)(-0.16,-0.16)
}
}
\put(2.5,2.5){$D_0:$}
\put(11.5,2.5){$D_1^2:$}
\end{picture}
\end{center}

\begin{center}
\setlength{\unitlength}{4mm}
\begin{picture}(20,5.8) 
\linethickness{0.6pt}
\put(5,-0.2){
\put(0,0){\OF}
\put(1,2){\GRF}
\put(0,4){\GR}
\qbezier(0,2)(0,3)(0,4)
\qbezier(2,6)(2,5)(2,4)
\qbezier(2,2)(2,1)(2,0)
\put(2,6){
\qbezier(0,0)(0.08,-0.08)(0.16,-0.16)
\qbezier(0,0)(-0.08,-0.08)(-0.16,-0.16)
}

}
\put(14,-0.2){
\put(0,0){\GRF}
\put(0,4){\GR}
\put(1,2){\GRF}
\qbezier(0,2)(0,3)(0,4)
\qbezier(2,6)(2,5)(2,4)
\qbezier(2,2)(2,1)(2,0)
\put(2,6){
\qbezier(0,0)(0.08,-0.08)(0.16,-0.16)
\qbezier(0,0)(-0.08,-0.08)(-0.16,-0.16)
}
}
\put(2.5,2.5){$D_1^1:$}
\put(11.5,2.5){$D_2:$}
\end{picture}
\end{center}
Denote by $\bar{D}_0$, $\bar{D}_1^1$, $\bar{D}_1^2$ and $\bar{D}_2$ the diagrams obtained from $\bar{D}$, after replacing the subdiagram $D$  by $D_0$, $D_1^1$, $D_1^2$ and $D_2$, respectively.
Then we have that the complex $C_n(\bar{D})$ associated to the diagram $\bar{D}$ is the total complex of the following complex of complexes:

\begin{equation*}
\begin{array}{ccccc}
&&C_n(\bar{D}_1^1)[1]\{-1\}&&\\
&\nearrow&&\searrow &\\
C_n(\bar{D}_0)&&\bigoplus&&C_n(\bar{D}_2)[2]\{-2\}\\
&\searrow&&\nearrow &\\
&&C_n(\bar{D}_1^2)[1]\{-1\}&&
\end{array}
\end{equation*}
Like previously, we have that 
\begin{equation}
C_n(\bar{D}_1^2)\cong C_n(\bar{D}_0)\{1\} \oplus C_n(\bar{D}_0)\{-1\},
\label{dirsum1}
\end{equation} 
and the projection of the map from $C_n(\bar{D}_0)$ onto the f\mbox{}irst summand of $C_n(\bar{D}_1^2)[1]\{-1\}$ is an isomorphism. Hence, again we can quotient by an acyclic complex and obtain that $C_n(\bar{D})$ is quasi-isomorphic to the total complex of:
\begin{equation}
\begin{array}{ccc}
C_n(\bar{D}_1^1)[1]\{-1\}&&\\
&\searrow&\\
\bigoplus&&C_n(\bar{D}_2)[2]\{-2\} \\
&\nearrow&\\
C_n(\bar{D}_0)[1]\{-2\}&&
\end{array}
\label{kompl1}
\end{equation}
Also, we have that 
\begin{equation}
C_n(\bar{D}_2)\cong C_n(\bar{D}_0) \oplus C_n(\bar{D}_3),
\label{dirsum2}
\end{equation}
(Proposition 33 in \cite{kovroz} -- ``categorif\mbox{}ication" of the last MOY) and that the map from the second summand of (\ref{dirsum1}) to the f\mbox{}irst summand of (\ref{dirsum2}) is an isomorphism (see \cite{kovroz} -- invariance under the Reidemeister III move). Here by $\bar{D}_3$ we denoted the diagram that is the same as $\bar{D}$ with the subdiagram $D$ replaced by the following diagram (see \cite{kovroz} and \cite{moy}):

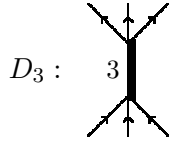
\begin{center}
\setlength{\unitlength}{3.5mm}
\begin{picture}(22,5) 
\linethickness{0.6pt}
\put(9,0)
{\qbezier(0,0)(0.75,0.75)(1.5,1.5)
\qbezier(1.5,0)(1.5,0.75)(1.5,1.5)
\qbezier(3,0)(2.25,0.75)(1.5,1.5)
\qbezier(0,5)(0.75,4.25)(1.5,3.5)
\qbezier(1.5,5)(1.5,4.25)(1.5,3.5)
\qbezier(3,5)(2.25,4.25)(1.5,3.5)
\put(0.7,2.2){\small{3}}
\linethickness{3pt}
\qbezier(1.5,1.5)(1.5,2.5)(1.5,3.5)
\linethickness{0.6pt}
\put(1.5,0.7){
\qbezier(0,0)(-0.08,-0.08)(-0.16,-0.16)
\qbezier(0,0)(0.08,-0.08)(0.16,-0.16)
}
\put(1.5,4.5){
\qbezier(0,0)(-0.08,-0.08)(-0.16,-0.16)
\qbezier(0,0)(0.08,-0.08)(0.16,-0.16)
}
\put(0.7,0.7)
{\qbezier(0,0)(0,-0.1)(0,-0.2)
\qbezier(0,0)(-0.1,0)(-0.2,0)
}
\put(2.5,4.5)
{\qbezier(0,0)(0,-0.1)(0,-0.2)
\qbezier(0,0)(-0.1,0)(-0.2,0)
}
\put(2.3,0.7)
{\qbezier(0,0)(0,-0.1)(0,-0.2)
\qbezier(0,0)(0.1,0)(0.2,0)
}
\put(0.5,4.5)
{\qbezier(0,0)(0,-0.1)(0,-0.2)
\qbezier(0,0)(0.1,0)(0.2,0)
}
\put(-3,2.2){$D_3:$}
}

\end{picture}
\end{center}

Hence $C_n(\bar{D})$ is quasi-isomorphic to the total complex of 
\begin{equation}
C_n(\bar{D}_1^1)[1]\{-1\} \to C_n(\bar{D}_3)[2]\{-2\}
\label{zezanje}
\end{equation}
On the other hand the complex $C_n(\bar{D}_3)$ is quasi-isomorphic to the total complex of both of the following two complexes:
$$C_n(\bar{D}_0)[-1]\to C_n(\bar{D}_2),$$
and
$$C_n(\bar{D}'_0)[-1]\to C_n(\bar{D}'_2),$$
where $\bar{D}'_0$ and $\bar{D}'_2$ are the diagrams obtained from $\bar{D}$ after replacing $D$ by following two diagrams, respectively:

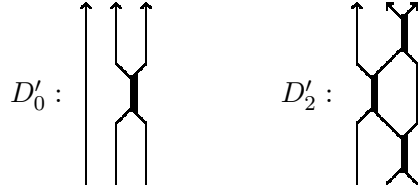
\begin{center}
\setlength{\unitlength}{4mm}
\begin{picture}(20,5.8) 
\linethickness{0.6pt}
\put(5,-0.2){
\put(1,0){\OF}
\put(1,4){\OF}
\put(1,2){\GRF}
\qbezier(0,6)(0,3)(0,0)
\put(0,6){
\qbezier(0,0)(0.08,-0.08)(0.16,-0.16)
\qbezier(0,0)(-0.08,-0.08)(-0.16,-0.16)
}
\put(1,6){
\qbezier(0,0)(0.08,-0.08)(0.16,-0.16)
\qbezier(0,0)(-0.08,-0.08)(-0.16,-0.16)
}
\put(2,6){
\qbezier(0,0)(0.08,-0.08)(0.16,-0.16)
\qbezier(0,0)(-0.08,-0.08)(-0.16,-0.16)
}
}
\put(2.5,2.5){$D'_0:$}
\put(14,-0.2){
\put(1,0){\GRF}
\put(1,4){\GR}
\put(0,2){\GRF}
\qbezier(2,2)(2,3)(2,4)
\qbezier(0,6)(0,5)(0,4)
\qbezier(0,2)(0,1)(0,0)
\put(0,6){
\qbezier(0,0)(0.08,-0.08)(0.16,-0.16)
\qbezier(0,0)(-0.08,-0.08)(-0.16,-0.16)
}
}
\put(11.5,2.5){$D'_2:$}
\end{picture}
\end{center}

Thus, the total complex of (\ref{kompl1}) (and hence $C_n(\bar{D})$) is quasi-isomorphic to the total complex of the following complex:

\begin{equation}
\begin{array}{ccc}
C_n(\bar{D}_1^1)[1]\{-1\}&&\\
&\searrow&\\
\bigoplus&&C_n(\bar{D}'_2)[2]\{-2\} \\
&\nearrow&\\
C_n(\bar{D}'_0)[1]\{-2\}&&
\end{array}
\label{kompl2}
\end{equation}

Thus, if we denote by $D_n$ the above complex shifted down in homological degree by 1, then we have that $C_n(\bar{D})\sim D_n[1]$, and all homology groups of $D_n$ are in nonnegative homological degrees, as required (the crossings in three diagrams appearing in (\ref{kompl2}) are all positive). We can now iterate the argument for each instance when an R2 move would occur in the $sl(2)$ case, since the two lower rightmost strands are both connected to the same thick edge in all three diagrams $D_1^1$, $D'_0$ and $D'_2$, and hence we can continue the process like with the initial diagram $D$.\\
\indent Completely analogously, we obtain the same result for the diagram $D$ of the form:

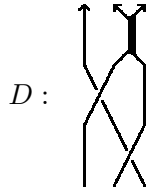
\begin{center}
\setlength{\unitlength}{4mm}
\begin{picture}(20,5.4) 
\linethickness{0.6pt}
\put(9,-0.2){
\put(1,0){\X}
\put(1,4){\GR}
\put(0,2){\X}
\qbezier(2,2)(2,3)(2,4)
\qbezier(0,6)(0,5)(0,4)
\qbezier(0,2)(0,1)(0,0)
\put(0,6){
\qbezier(0,0)(0.08,-0.08)(0.16,-0.16)
\qbezier(0,0)(-0.08,-0.08)(-0.16,-0.16)
}
}
\put(6.5,2.5){$D:$}
\end{picture}
\end{center}

Thus, we obtain the required shift in homological degree. \kraj

\end{document}